\documentclass[12pt]{article}
\usepackage{graphicx}
\usepackage{amsmath}
\usepackage{amssymb}

\usepackage{cite}

\def\vr{\varrho}
\def\el{\ell}

\def\nn{\nonumber}
\def\a{\alpha}

\def\s{\sigma}
\def\la{\lambda}

\def\vp{\varphi}

\def\ve{\varepsilon}
\def\wh{\widehat}
\def\wt{\widetilde}
\def\ov{\overline}

\def\p{\partial}

\def\BC{{\mathbb C}}
\def\BR{{\mathbb R}}
\def\BN{{\mathbb N}}

\def\clb{{\mathcal B}}

\def\diag{\mathrm{diag}}
\newcommand{\E}{\mathrm{e}}
\newcommand{\I}{\mathrm{i}}

\newtheorem{Pa}{Paper}[section]

\newtheorem{Rk}[Pa]{{\bf Remark}}

\newtheorem{Ee}[Pa]{{\bf Example}}

\newtheorem{Pn}[Pa]{{\bf Proposition}}

\title{Generalized B\"acklund-Darboux transformation: conservation laws, \\ rational extensions and bispectrality}

\author{Alexander Sakhnovich}

\date{}

\begin{document}
\maketitle

\begin{abstract} B\"acklund-Darboux transformations are closely related to the integrability
and symmetry problems. For the generalized B\"acklund-Darboux transformation (GBDT),
we consider  conservation laws, rational extensions and bispectrality. We use the case of the
nonlinear optics equation (and its auxiliary linear system) as an example.
\end{abstract}

\section{Introduction}
\setcounter{equation}{0}
 B\"acklund-Darboux transformations  and related commutation methods constitute one of the most  fruitful 
 approaches to the construction of explicit solutions of linear and nonlinear equations and to the
 explicit spectral theoretic results 
 (see, e.g., \cite{Ci, Ge, GeT, Gu, KoSaTe, MS} and numerous references therein).
  Similar to various symmetries studied in group analysis,
 B\"acklund-Darboux transformations (BDTs) transform solutions 
 of  linear differential equations of some fixed class
 into solutions of differential equations of the same (or another fixed) class and potentials into potentials, 
 whereas potentials are often solutions of the
 corresponding integrable nonlinear equations.
 Some connections between Darboux transformations and potential symmetries are
 discussed in an interesting paper \cite{PKI}.
 
  Though the idea of BDTs and commutation methods goes back to 
 B\"ack-lund, Darboux and Jacobi,
 the notion of the so called {\it Darboux matrix} has a much later origin. 
 According to J.~Cie\'sli\'nski   \cite{Ci},
"all approaches to the construction of Darboux matrices originate in the dressing method".
 In order to explain the notion of the Darboux matrix, we consider some initial linear system 
 \begin{equation}\label{1}
u_x=G(x,z) u, \quad u\in \BC^m, 
\end{equation}
where $u_x:=\frac{d}{dx}u\,$ ($u_x=\frac{\p}{\p x}u$ in the case of several variables), $G(x,z) $ is some 
$m \times m$ matrix function (i.e., $G(x,z) \in \BC^{m\times m}$) and $\BC$ stands for
the complex plain.
 The \textit{Darboux matrix} (matrix function) $\wh w$ satisfies the equation
 \begin{equation}\label{2}
\wh w_x=\wt G \wh w-\wh wG , \quad \wt G(x,z) \in \BC^{m\times m}.
\end{equation}
It is easy to see that if \eqref{1} and \eqref{2} hold,  the product 
$$\wt u=\wh wu$$ 
satisfies the  {\it transformed system}
 $$\wt u_x=\wt G(x, z)\wt u.$$
The cases, where the Darboux matrix $\wh w$ can be constructed explicitly, if $u$ and $G$ are given
explicitly, are of special interest.
Clearly, in those cases we construct explicitly $\wt u$ and also $\wt G$.
 
 We are interested in  a generalized version of the B\"acklund-Darboux transformation for which we use the acronym GBDT.
 In GBDT (for the case of one space variable) Darboux matrix is presented as the transfer matrix function
 in Lev Sakhnovich form:
  \begin{align} \label{3} &
\wh w(z)=w_A(z)=  I_m-\Pi_2^*S^{-1}(A_1-z I_n)^{-1}\Pi_1; \quad \Pi_k \in \BC^{n \times m}, \quad S\in \BC^{n \times n};
\\ & \label{4} {\mathrm{where}} \quad \quad
A_1S-SA_2=\Pi_1\Pi_2^*,
\end{align}
$I_m$ is the $m\times m$ identity matrix. We note that an additional variable $x$ appears in $ w_A$, $\Pi_1$, $\Pi_2$ and $S$,
so that we deal with matrix functions $\Pi_k(x)$, $S(x)$ and Darboux matrix $w_A(x,z)$.
GBDT is applicable to all systems
$u_x=G(x,z)u$, where  $G$ rationally depends on $z$, that is,
$$
G(x,z)=\Big(\sum_{k=0}^r z^k q_k(x)+\sum_{s=1}^j \sum_{k=1}^{r_s}(z -c_s)^{-k}
q_{sk}(x)\Big)
$$
(see the paper \cite{ALS1}, the reviews \cite{ALS2, ALS2'} and Chapter 7 in \cite{SaSaR} as well as numerous
references therein).
The coefficients $\wt q_k$ and $\wt q_{sk}$ of the transformed (rational) matrix function $\wt G(x,z)$ are  expressed via the
coefficients $q_k$ and $ q_{sk}$ of  $ G(x,z)$ and matrix functions $S(x)$, $\Pi_1(x)$ and $\Pi_2(x)$.
Matrix functions $\Pi_k(x)$ are constructed as solutions of systems related
to the initial system.

Here, we shall consider in detail GBDT for the system
\begin{align}\label{5} &
u_x=G(x,z)u, \quad G=\I z D-[D,\vr], \quad D={\mathrm{diag}}\, \{d_1,d_2, \ldots,d_m\}=D^*;
\\  \label{6} &
\vr^*=B  \vr B, \quad B={\mathrm{diag}}\, \{b_1,b_2, \ldots,b_m\} \quad (b_k=\pm 1); \quad [D,\vr]:=D\vr -\vr D,
\end{align}
where diag denotes a diagonal matrix. System \eqref{5} is an auxiliary system for the well-known  $N$-wave (nonlinear optics)
equation:
\begin{align} \label{7} &    
[D,  \vr_t ] - [\wh D,  \vr_x ]
=
\big[[D, \vr],\,[\wh D, \vr]\big],   
\\ &     \label{8} 
\wh D={\mathrm{diag}}\, \{\wh d_1, \ldots , \wh d_m\}=\wh D^*.
\end{align}
More precisely, equation \eqref{7} is equivalent to the compatibility condition
 \begin{align} \label{9} &
G_t-F_x+[G,F]=0
\end{align}
of the auxiliary systems
\begin{align}\label{10} &
u_x=G(x,t,z)u, \quad G=\I z D-[D,\vr], 
\\  \label{11} &
u_t=F(x,t, z)u, \quad F=\I z \wh D-[\wh D,\vr].
\end{align}
We shall  study conservation laws, rational extensions  and bispectrality for the case
of GBDT for system \eqref{5} and also for the $N$-wave equation \eqref{7}.

The integrable model describing interaction of  three wave packages \cite{ZaMa}  is the most well-known subcase of the $N$-wave equation \eqref{7}.
Putting, for instance, $m=3$, $B=I_3$ and $d_1>d_2>d_3$, and using transformations \cite[Ch. 3]{NoMa}
\begin{align}\nn
\psi_1=(\wh d_2 -\wh d_1)/(d_1 - d_2), \quad & \psi_2=(\wh d_3 -\wh d_2)/(d_2 - d_3), \\
 \label{R1} &     \psi_3=(\wh d_3 -\wh d_1)/(d_1 - d_3);
\\  \nn
\phi_1=-\I \sqrt{d_1-d_2}\, \vr_{12}, \quad & \phi_2=-\I \sqrt{d_2-d_3}\, \vr_{23}, \\
&     \label{R2} 
 \phi_3=-\I \sqrt{d_1-d_3}\, \vr_{13},
\end{align}
we rewrite \eqref{7} in the standard form of the corresponding 3-wave interaction:
\begin{align} \nn  
(\phi_1)_t+\psi_1 (\phi_1)_x=\I \ve \ov \phi_2 \phi_3, \quad & (\phi_2)_t+\psi_2 (\phi_2)_x=\I \ve \ov \phi_1 \phi_3, \\
\label{R3} &   (\phi_3)_t+\psi_3 (\phi_3)_x=\I \ve  \phi_1 \phi_2,
\end{align}
where $\ov \phi(x,t)$ stands for the function, which takes values complex conjugate to $\phi(x,t)$, and
\begin{align} \nn
\ve=&(d_1\wh d_2-d_2\wh d_1+d_2\wh d_3-d_3\wh d_2 + d_3\wh d_1-d_1\wh d_3)
\\ \label{R4} & \quad \times    
\big((d_1-d_2)(d_1-d_3)(d_2-d_3)\big)^{-1/2}.
\end{align}
 If $\a$ is a scalar value or matrix, the notation $\ov \a$ stands for the scalar, which is complex conjugate to $\a$,
 or the matrix with entries complex conjugate to the entries of $\a$, respectively.
\section{Preliminaries}
\setcounter{equation}{0}
GBDT for the system \eqref{5} and for the $N$-wave equation was described
in \cite{ALS1} (see also \cite[Subsection 1.1.3 and Section 7.1]{SaSaR} and references therein).
In this section, we give some necessary definitions and results from \cite{ALS1} and  \cite[Subsection 1.1.3 and Section 7.1]{SaSaR}.
We consider the case $x\geq 0$ for the system \eqref{5} and the case $x \geq 0, \, t \geq 0$ for the $N$-wave equation. 
We fix $n \times n$ matrices $A$ and $S(0)=S(0)^*$
and  an $n \times m$ matrix $\Pi(0)$ such that
\begin{equation}       \label{12}
AS(0)-S(0)A^*=\I\Pi(0)B \Pi(0)^*. 
\end{equation}
We introduce  the $n \times m$ matrix function 
$\Pi(x)$ and the $n \times n$ matrix function $S(x)$ via initial values $\Pi(0)$ and $S(0)$ and differential equations
\begin{align}      &\label{13}
\Pi_x=-\I A\Pi D+\Pi[D, \vr], \quad S_x=\Pi D B \Pi^*.
\end{align}
Relations \eqref{12} and \eqref{13} yield $AS(x)-S(x)A^*=\I\Pi(x)B \Pi(x)^*$.
Then the following proposition describes GBDT determined by $A$, $S(0)$
and   $\Pi(0)$.
\begin{Pn}\label{PnGTNW}
Let  system \eqref{5} $($such that  \eqref{6} holds$)$ be given. Then $($in the points of invertibility
of $S(x)$$)$ the matrix function
\begin{equation}       \label{14}
w_A(x, z)=  I_m-\I B \Pi(x)^*S(x)^{-1}(A-zI_n)^{-1}\Pi(x)
\end{equation}
is a Darboux matrix of \eqref{5} and satisfies the equation
\begin{equation}       \label{15}
\frac{d}{dx}w_A(x, z)=\big(\I zD-[D, \wt \vr(x)]\big)w_A(x, z)-
w_A(x, z)\big(\I zD-[D, \vr(x)]\big),
\end{equation}
where
\begin{equation}       \label{16}
\wt \vr=\vr  -B\Pi^*S^{-1}\Pi, \quad \wt \vr^*=B \wt \vr B.
\end{equation}
Moreover, if $\det S(x)\not=0$ for $x \geq 0$, a normalized fundamental solution $\wt w(x,z)$ of the transformed system
\begin{align}\label{17} &
\wt w_x=(\I z D-[D, \wt \vr])\wt w,  \quad \wt w(0,z)=I_m
\end{align} 
is given by the equality
\begin{align}\label{18} &
\wt w(x,z)=w_A(x,z)w(x,z) w_A(0,z)^{-1},
\end{align} 
where $w$ is the fundamental solution of the initial system \eqref{5} normalized by $w(0,z)=I_m$.
\end{Pn}
Next, we  add the variable $t$ and consider $\Pi(x,t)$, $S(x,t)$ and $w_A(x,t,z)$ which are determined by $A$, $\Pi(0,0)$ and $S(0,0)=S(0,0)^*$
via equations \eqref{13}, \eqref{14} and
\begin{align}      &\label{19}
\Pi_t=-\I A\Pi \wh D+\Pi[\wh D, \vr], \quad S_t=\Pi \wh D B \Pi^*.
\end{align}
Instead of \eqref{12} we assume that
\begin{equation}       \label{20}
AS(0,0)-S(0,0)A^*=\I\Pi(0,0)B \Pi(0,0)^*,
\end{equation}
which implies that  $AS(x,t)-S(x,t)A^*=\I\Pi(x,t)B \Pi(x,t)^*$. Proposition \ref{PnGTNW} yields:
\begin{Pn}\label{PnNW} Let an $m \times m$ matrix function $\vr$ $(\vr^*=B  \vr B)$ be continuously
differentiable and satisfy the $N$-wave equation \eqref{7}. Then $\wt \vr$ of the form
\begin{equation}       \label{21}
\wt \vr(x,t):=\vr(x,t)  -B\Pi(x,t)^*S(x,t)^{-1}\Pi(x,t)
\end{equation}
satisfies $($in the points of invertibility of $S)$ the equality $\wt \vr^*=B  \wt \vr B$ and the $N$-wave equation.
\end{Pn}
Using \cite[Therem 6.1]{SaSaR} on wave functions, we easily obtain the next statement.
\begin{Rk} \label{Rk1} 
If the conditions of Proposition \ref{PnNW} hold,
 $\det S(x,t)\not=0$ on the semi-band $0 \leq x <\infty, \, \,  0\leq t \leq \ve$ and $w(x,t,z)$ is the initial wave
function $($i.e., $w_x=Gw$, $w_t=Fw$ and  $w(0,0,z)=I_m)$, then the transformed wave function $\wt w(x,t,z)$ is given by the
equality
\begin{equation}       \label{22}
\wt w(x,t,z)=w_A(x,t,z)w(x,t,z)w_A(0,0,z)^{-1}.
\end{equation}
\end{Rk}
\begin{Ee} \label{REe}
The real-valued case is of special interest \cite[Section 3.4]{NoMa}. It is immediate from Proposition \ref{PnNW}
$($and formulas defining $\wt \vr$ considered in this proposition$)$ that equalities
\begin{equation}       \label{R5}
\vr = \ov \vr, \quad  A = - \ov A, \quad  \Pi(0,0)=\ov \Pi(0,0), \quad  S(0,0)=\ov S(0,0),
\end{equation}
yield the equality $\displaystyle{\wt \vr = \ov{\wt \vr}}$. Assume that \eqref{R5} holds, and so $\displaystyle{\wt \vr = \ov{\wt \vr}}$. 
Then, in the subcase $m=3$, $B=I_3$, $d_1>d_2>d_3$ of the
3-wave equation, we can rewrite \eqref{R3} in the real-valued form and  obtain solutions of the corresponding exact resonance
equations. Namely, we set
\begin{align}  &     \nn
\vp_1=-\I \phi_1=- \sqrt{d_1-d_2}\, \wt \vr_{12}, \quad \vp_2=-\I \phi_2=- \sqrt{d_2-d_3}\, \wt \vr_{23}, \\
\nn &  \vp_3=-\I \phi_3=- \sqrt{d_1-d_3}\, \wt \vr_{13},
\end{align}
and, taking into account \eqref{R3}, derive
\begin{align} \nn
(\vp_1)_t+\psi_1 (\vp_1)_x= \ve \vp_2 \vp_3, \quad & (\vp_2)_t+\psi_2 (\vp_2)_x= \ve  \vp_1 \vp_3, \\
\label{R6} &    
 (\vp_3)_t+\psi_3 (\vp_3)_x=- \ve  \vp_1 \vp_2.
\end{align}
\end{Ee}
\section{Conservation laws, rational extensions  and bispectrality}
\setcounter{equation}{0}
\subsection{Conservation laws}
When $\vr$ is differentiable with respect to $t$ (to $x$), we can differentiate both sides of the first
relation in \eqref{13} (in \eqref{19}) with respect to $t$ (to $x$). If $\vr$ is differentiable with respect to $x$ and $t$ and $\vr_t$ is continuous, then, 
according to the so called Clairaut's (or Schwarz's) theorem, we have $\Pi_{xt}=\Pi_{tx}$. Thus $\Pi_{xt}=\Pi_{tx}$ is the necessary condition of
the compatibility of the first relations in \eqref{13} and \eqref{19}. In a similar way, $S_{xt}=S_{tx}$ is a necessary condition of the
compatibility of  the second relations in \eqref{13} and \eqref{19}. In other words, the necessary compatibility conditions are
\begin{align}\label{c1} &
(-\I A\Pi D+\Pi[D, \vr])_t=(-\I A\Pi \wh D+\Pi[\wh D, \vr])_x,
\\ \label{c2} &
(\Pi  D B \Pi^*)_t=(\Pi \wh D B \Pi^*)_x,
\end{align} 
where $\Pi_x$ should be substituted by $-\I A\Pi D+\Pi[D, \vr]$ and $\Pi_t$ should  be substituted by $-\I A\Pi \wh D+\Pi[\wh D, \vr]$.
\begin{Pn} \label{NesCond} Let $\vr$ be differentiable and satisfy the $N$-wave equation \eqref{7}. Then the equalities \eqref{c1} 
and \eqref{c2} hold.
\end{Pn}
{\bf Proof.} Substituting expressions for $\Pi_t$  and $\Pi_x$ we obtain 
\begin{align}\nonumber 
(-\I A\Pi D+\Pi[D, \vr])_t= & -\I A(-\I A\Pi \wh D+\Pi[\wh D, \vr])D
\\ \nn &
+(-\I A\Pi \wh D+\Pi[\wh D, \vr])[D, \vr]+\Pi[D, \vr_t],
\\ \nonumber 
(-\I A\Pi \wh D+\Pi[\wh D, \vr])_x=& -\I A(-\I A\Pi  D+\Pi[ D, \vr])\wh D
\\ \nn &
+(-\I A\Pi  D+\Pi[ D, \vr])[\wh D, \vr]+\Pi[\wh D, \vr_x].
\end{align} 
Therefore, taking into account \eqref{7} we derive
\begin{align}\nonumber &
(-\I A\Pi D+\Pi[D, \vr])_t-(-\I A\Pi \wh D+\Pi[\wh D, \vr])_x
\\ \nn & \quad
=\Pi\big([D,  \vr_t ] - [\wh D,  \vr_x ]
-
\big[[D, \vr],\,[\wh D, \vr]\big]\big)=0,
\end{align} 
and so \eqref{c1} is proved. Using the same substitutions as before, we have also 
\begin{align}\label{c3} &
(\Pi  D B \Pi^*)_t= \big(-\I A\Pi \wh D+\Pi[\wh D, \vr]\big)DB\Pi^*+\Pi DB\big(\I \wh D \Pi^* A^*+[\vr^*,\wh D]\Pi^*\big),
\\ \label{c4} &
(\Pi \wh D B \Pi^*)_x=\big(-\I A\Pi  D+\Pi[ D, \vr]\big)\wh DB\Pi^*+\Pi \wh DB\big(\I  D \Pi^* A^*+[\vr^*, D]\Pi^*\big).
\end{align} 
Equality \eqref{c2} follows from \eqref{c3}, \eqref{c4} and $B\vr^*=\vr B$.
$\Box$

The sufficiency of the proved above equalities $\Pi_{xt}=\Pi_{tx}$ and $S_{xt}=S_{tx}$ for the compatibility
of \eqref{13} and \eqref{19} is not self-evident but may be proved in a way similar to the proof of \cite[Theorem 6.1]{SaSaR}.
\begin{Rk}
We note that relations \eqref{c1} and \eqref{c2} may be considered as conservation laws for $\Pi(x,t)$.
\end{Rk}
\subsection{Rational extensions}
According to Propositions \ref{PnGTNW} and \ref{PnNW}, in the case 
\begin{align}      &\label{23}
\vr \equiv 0, \quad \sigma(A)=\{0\} \quad (\sigma \, {\mathrm{stands \, for \, spectrum}}),
\end{align}
the transformed potentials $\wt \vr$ are rational matrix functions ({\it rational extensions} in the terminology
of \cite{Gran}). More precisely, if $\vr \equiv 0$ and $A$ is nilpotent, formulas \eqref{13}, \eqref{14}, \eqref{16}, \eqref{18}, \eqref{19}, \eqref{21} and 
\eqref{22} imply
 the following proposition.
 \begin{Pn}\label{PnAlgExt} Let \eqref{23} hold. Then $\Pi(x)$ and $S(x)$ are matrix polynomials with respect to $x$,
$\Pi(x,t)$ and $S(x,t)$ are matrix polynomials with respect to $x$ and $t$; $\wt \vr(x)=p(x)/\det S(x)$ and $\wt \vr(x,t)=p(x,t)/\det S(x,t)$,
where $p(x)$ is a matrix polynomial with respect to $x$, $p(x,t)$ is a matrix polynomial with respect to $x$ and $t$,
$\det S(x)$ is a polynomial with respect to $x$ and $\det S(x,t)$ is a polynomial with respect to $x$ and $t$.
Moreover, for Darboux matrices $w_A$ we have $w_A(x,z)=P(x,z)/\det S(x)$ and $w_A(x,t,z)=P(x,t,z)/\det S(x,t)$,
where $P(x,z)$ $\,(P(x,t,z))$ is a matrix polynomial with respect to $x$ and $1/z$ $\, (x,\, t$ and $1/z)$. Thus,
the transformed fundamental solutions $\wt w(x,z)$ and wave functions $\wt w(x,t, z)$ are expressed via matrix polynomials and
exponents of diagonal matrices:
 \begin{align}    \nn
\wt w(x,z) &=w_A(x,z)\E^{\I z x D} w_A(0,z)^{-1}, 
\\    \label{24} 
\wt w(x,t,z)&=w_A(x,t, z)\E^{\I z( x D+t \wh D)} w_A(0,0,z)^{-1}.
\end{align}
\end{Pn}
{\bf Proof.} Since $\vr \equiv 0$, the first equation in \eqref{13} yields that the $i$th column $f_i$ of $\Pi$ is given
by $f_i(x)=\E^{-\I d_i x A}f_i(0)=\sum_{k=0}^{n-1}\frac{1}{k!}(-\I d_i x A)^kf_i(0)$, and so $\Pi(x)$ is a matrix polynomial.
(Here we used the equality $A^n=0$, which holds for all $n \times n$ nilpotent matrices.) Hence, the second
equation in \eqref{13} yields that $S(x)$ is a matrix polynomial. In the same way we prove the polynomial form of
$\Pi(x,t)$ and $S(x,t)$, and the required properties of $\wt \vr$ follow. 
According to formula (1.84)
from \cite[p. 24]{SaSaR} we have $w_A(z)Bw_A(\ov z)^*=B$ and, in particular,
the equalities
$$w_A(0,z)^{-1}=Bw_A(0,\ov z)^*B, \quad w_A(0,0,z)^{-1}=Bw_A(0,0,\ov z)^*B$$
are valid.
Finally, the representations of
$w_A(x,z)$ and $w_A(x,t,z)$ follow from the definition \eqref{14}, properties of $\Pi$ and $S$ and the
expansion $(A-z I_n)^{-1}=-z^{-1}\sum_{k=0}^{n-1}(z^{-1}A)^k$. 
$\Box$
\begin{Rk} System \eqref{5} includes important subclasses:  self-adjoint and skew-self-adjoint Dirac-type systems
$($which are also called Zakharov-Shabat or AKNS systems$)$. Namely,  putting
 \begin{equation}       \label{25}
D=j=\begin{bmatrix}I_{m_1} & 0 \\  0 & -I_{m_2}\end{bmatrix} \quad (m_1+m_2=m), \quad V=\begin{bmatrix}0 & v \\  v^* & 0\end{bmatrix}
\end{equation}
and $B=j$, we rewrite  \eqref{5}, \eqref{6} in the form of the  self-adjoint Dirac-type system
 \begin{equation}       \label{26}
u_x=\I(zj+jV)u, \quad v=2\I\vr_{12},
\end{equation}
where $\vr_{12}$ is the upper right  $m_1 \times m_2$ block of $\vr$. Assuming \eqref{25} and $B=I_m$,
we rewrite  \eqref{5}, \eqref{6} in the form of the  skew-self-adjoint Dirac-type system
 \begin{equation}       \label{27}
u_x=(\I zj+jV)u, \quad v=-2\vr_{12}.
\end{equation}
\end{Rk}
\subsection{Bispectrality}
The notion of bispectrality was introduced in \cite{DG} (see also \cite{Grun} and references therein).
Matrix bispectrality for system \eqref{5} was studied in \cite{SZ}. Namely, it was shown that for $\wt \vr$
constructed in  Proposition \ref{PnAlgExt} and, correspondingly, for the solution $W(x,z)= w_A(x,z) \E^{\I z x D}$
of system 
 \begin{equation}       \label{b0}
W_x=(\I z D -[D,\wt \vr])W,
\end{equation}
there is a non-degenerate matrix linear differential operator $B$ with respect to the variable $z$
such that $B(z)W=0$. It was suggested in \cite{Grun} (and in some earlier works by F. Gr\"unbaum and coauthors) that the bispectrality
requirement $W\clb(z)=\Theta(x)W$ is more meaningful than the requirement $B(z)W=\Theta(x)W$.
Unfortunately, the approach from \cite{SZ} does not work properly for that case. Therefore, we should consider the
system
 \begin{equation}       \label{b1}
W_x=(\I z D -[D,\wt \vr])W- \I z W D
\end{equation}
instead of system \eqref{5}. In view of \eqref{15}, the matrix function $w_A$ satisfies this system in the case $\vr \equiv 0$.
Thus, we can deal in the same way as in \cite{GKS2}, where Dirac systems were treated. Namely, if $\s(A)$
is concentrated at one point $\la\in \BC$
we write down the resolvent $(A-z I_n)^{-1}$
in the form
\begin{equation}       \label{b2}
(A-z I_n)^{-1}=-\sum_{k=1}^n(z-\la)^{-k}(A-\la I_n)^{k-1}.
\end{equation}
By virtue of \eqref{14} and \eqref{b2} the next proposition is immediate.
\begin{Pn} Assume that $\vr \equiv 0$, 
the spectrum of $A$ is concentrated
at some point $\la\in \BC$ and $\wt \vr$ is given by \eqref{16}.
Then system \eqref{b1} is bispectral in the sense of   \cite{Grun}.
\end{Pn}
Indeed, it is easy to choose some coefficients $c_s$ so that the operator
$B(z)=\sum_{s=1}^{n+1} c_s(z- \la)^s\frac{\p^s}{\p z^s}$ is non-degenerate and satisfies the equalities
\begin{equation}       \label{R7}
B(z) (z- \la)^{-k}=0 \qquad (1\leq k \leq n).
\end{equation}
According to formula \eqref{15} and identity $\vr \equiv 0$, the matrix function $W=w_A(x,z)$ satisfies \eqref{b1}.
Therefore, definition \eqref{14} of $w_A$ and equalities \eqref{b2} and \eqref{R7} imply $B(z)W=0$.

If we want to consider system \eqref{b0}, where $\vr \equiv 0$ and $\s(A)=\{\la\}, \quad \la \in \BR$,
we rewrite the first relation in \eqref{13} in the form
 \begin{equation}       \label{29}
\left(\Pi(x)\E^{\I \la x D}\right)_x=-\I(A-\la I_m)\left(\Pi(x)\E^{\I \la x D}\right)D.
\end{equation}
Hence, taking into account the fact that  $\s(A)$ is concentrated at  $\la$, we see that  $\Pi(x)\E^{\I \la x D}=P(x,\la)=P(x)$,
where $P(x)$ is a matrix polynomial with respect to $x$. Substituting $\Pi(x)=P(x)\E^{-\I \la x D}$
into \eqref{13} and using the equalities $D=D^*, \, \la=\ov{\la}$, we obtain
 \begin{equation}       \label{30}
S_x=P(x)\E^{-\I \la x D}DB\E^{\I \la x D}P(x)^*=P(x)DBP(x)^*, \quad {\mathrm{i.e.,}} \quad S(x)=\wt P(x),
\end{equation}
where $\wt P(x)$ a matrix polynomial with respect to $x$.  Introduce an operator
$\clb=\diag\{\clb_1,\ldots , \clb_m\}$  such that  $\clb_k$ has the form
 \begin{equation}       \label{28}
f \,  \clb_k =\sum_{\el=1}^N c_{\el k}(z)\frac{\p^{\el}}{\p z^{\el}}\big((z- \la)^m f\big), \quad N\in \BN,
\end{equation}
and notice that the equality
\begin{align}  &     \label{31}
\sum_{\el=1}^N c_{\el k}(z)\frac{\p^{\el}}{\p z^{\el}}\left((z- \la)^m w_A(x,z) \E^{\I z x D}e_k\right)=0, \\
\nn & \qquad \qquad \qquad 
 e_k:=\{\delta_{ik}\}_{i=1}^m \quad (1\leq k \leq m)
\end{align}
 is equivalent  to the relation
 \begin{align}       \nonumber &
\sum_{\el=1}^N c_{\el k}(z)\left(\frac{\p}{\p z}+\I d_k x\right)^{\el}
\\ \nn & \quad \times
\left((z- \la)^m \E^{\I \la x D}\left(I_m-\I BP(x)^*\wt P(x)^{-1}(A-z I_n)^{-1} P(x)\right)\E^{-\I \la x D}e_k\right)=0,
\end{align}
which, in turn, is equivalent to
 \begin{align}     \label{32}  &
 \sum_{\el=1}^N c_{\el k}(z)\left(\frac{\p}{\p z}+\I d_k x\right)^{\el}
\\ \nn & \quad \times 
 \left((\det \wt P(x))(z- \la)^m \left(I_m-\I BP(x)^*\wt P(x)^{-1}(A-z I_n)^{-1} P(x)\right)e_k\right)=0.
\end{align}
{\bf Open Problem.} Are there some cases, where matrices $D$, $A$, $S(0)$ and $\Pi(0)$ generate nontrivial $\wt \vr$ such that
$W\clb(z)=0$ for the solution $W(x,z)= w_A(x,z) \E^{\I z x D}$ of \eqref{b0} and for $\clb\not=0$ of the form given above?
In other words, are there $D$, $A$, $S(0)$, $\Pi(0)$ and $\{c_{lk}\}$ satisfying \eqref{32}? 

{\bf Acknowledgement.}
The research was supported  by the  Austrian Science Fund (FWF) under Grant  No. P24301.

\vspace{0.5em}

Alexander Sakhnovich, University of Vienna, \\ oleksandr.sakhnovych@univie.ac.at

\end{document}